\documentclass[11pt]{article}
\usepackage{authblk}
\usepackage{graphicx}    % needed for including graphics e.g. EPS, PS
\usepackage{color}
\usepackage{amsmath,
            amsthm,
            amssymb}     % for mathematics
\usepackage[all]{xy}
\usepackage{epstopdf}
\usepackage{graphicx}
\newcommand{\beq}{\begin{equation}}
\newcommand{\eeq}{\end{equation}}
\newtheorem{thm}{Theorem}[section]
\newtheorem{cor}[thm]{Corollary}
\newtheorem{lem}[thm]{Lemma}
\newtheorem{prop}[thm]{Proposition}
\newenvironment{df}[1][Definition]{\begin{trivlist}
\item[\hskip \labelsep {\bfseries #1}]}{\end{trivlist}}
\newenvironment{ex}[1][Example]{\begin{trivlist}
\item[\hskip \labelsep {\bfseries #1}]}{\end{trivlist}}

\newcommand{\bpf}{\begin{proof}}
\newcommand{\epf}{\end{proof}}

\begin{document}
\title{The Classification of Two Dimensional topological Field Theories}
\author{Geoffrey Lee}
\affil{541 Leverett Mail Center, 28 Dewolf Street, Cambridge, MA 02138, USA.}
\affil{Email: glee@college.harvard.edu}
\maketitle
\begin{abstract}
The goal of this thesis is to introduce some of the major ideas behind extended topological quantum field theories with an emphasis on explicit examples and calculations. The statement of the Cobordism Hypothesis is explained and immediately used to classify framed and oriented extended $2$ dimensional topological quantum field theories. The passage from framed theories to oriented theories is equivalent to giving homotopy fixed points of an $SO(n)$ action on the space of field theories. This thesis then constructs extended $2$ dimensional Dijkgraaf-Witten theory (also called finite gauge theory) as an example of a $2$ dimensional extended field theory by assigning invariants at the level of points and extending up. Finally, it is concluded that Dijkgraaf-Witten theory is the only example of an extended framed $2$ dimensional topological quantum field theory by showing that any field theory is equivalent to Dijkgraaf-Witen theory for some cyclic group.
\end{abstract}
\clearpage
\tableofcontents
\clearpage
\setcounter{section}{-1}
\section{Introduction}
\subsection{Motivation}
In the last few decades, there has been considerable research interest in the subject of topological quantum field theories. After Atiyah's axiomitization of these topological field theories in his 1989 paper ``Topological quantum field theories'' \cite{atiyah88}, there has been much interest in classifying field theories, developing new ones calculating interesting manifold invariants, and finding physical models of such field theories. From a pratical perspective, topological quantum field theories describe systems capable of performing quantum computation. To date, the best candidate system to be described by a topological quantum field theory is the $\nu=5/2$ fractional quantum hall effect - a $3$ dimensional field theory \cite{moore91}. Heuristically the quasiparticles of this system are described by tiny tubes representing the movement of the particles through time; one would use this system to perform quantum computations by studying the braiding of these tubes in a system. The quantum amplitude of that time evolution would depend only on the regular isotopy type of the link; in other words, a small perturbation in the path of the particle will not change the quantum amplitude as long as the perturbation does not affect the way the tubes cross. In this system, the different isotopy types of links would correspond to different logic gates in computation. The benefit to such a ``topological'' model of quantum computing is that computation is not affected by the inevitable small amounts of ``noise'' leaking into all real-world systems. Because the paths of the particles only matter up to isotopy, a local perturbation of the path will not change the way these paths braid meaning that as long as the particles were ``far enough apart'' to begin with, the perturbed paths will still compute the same logic gate. Therefore, the real practical benefit to topological quantum computing is the ability to compute even when the system is not perfectly insulated \cite{nayak08}.

From a more mathematical perspective, topological field theories are interesting because they are tools to systematically organize and produce topological invariants. Already, $3$ dimensional field theories have been used to organize knot invariants and polynomials \cite{witten89}. The subject of this thesis however, will mostly be $2$ dimensional field theories: in particular, we produce a complete classification of all $2$ dimensional field theories by calculating the necessary and sufficient pieces of data required to define a $2$ dimensional field theory. We will see that these field theories nicely translate geometric conditions to algebraic conditions with the result that $2$ dimensional field theories are essentially special kinds of algebras over a field.

These geometric conditions, of course, are inspired by the physics of topological field theories. The spatial dimension is typically described by $M$, an $n$ dimensional manifold without boundary with bordisms of $n$-manifolds describing the evolution of this system. If there is a bordism $\emptyset\rightarrow M$ then this bordism selects the vacuum state of $M$. If $\partial M=\emptyset$, then $Z(M)$ is the vacuum expectation value \cite{dijkgraaf90}. In general, one can ask how the quantum system evolves over a period of time $t\in I$. An important aspect of topological field theories is that the quantum amplitude for time evolution is invariant under diffeomorphism - therefore, we can take the cylinder $M\times I$ to represent the ``identity map.'' On the other hand, the quantum nature of the system allows for ``quantum fusion'' and ``quantum splitting,'' the amplitudes for which are described by bordisms $M\rightarrow N$. All of this data can be reorganized in a coherent manner using the language of categories and functors:
\begin{df} An ``$n+1$ dimensional topological quantum field theory'' (or ``topological field theory'') is a symmetric monoidal functor $Z:\text{Cob}(n)\rightarrow \text{Vect}(k)$. \end{df}
This formalism is precisely what Aityah described in \cite{atiyah88}. In particular, the class of topological quantum field theories forms a category itself with morphisms being the natural transformations. It is sensible then, to ask if this category is equivalent to another category, or at least if there is a way to produce a list of distinct equivalence classes of field theories. This sort of classification is precisely what this thesis achieves in dimension $2$.

\subsection{Monoidal Categories}
The basic language with which to formulate topological field theories is that of monoidal categories and functors. We will give a brief overview of the relevant definitions and concepts in this section leaving the details to MacLane \cite{maclane98}.
\begin{df} A ``monoid'' is a set $S$ equipped with an associative operation $\circ:S\times S\rightarrow S$ such that there is a distinguished element $1$ satisfying the property that $s\circ 1 = 1 \circ s$ for all $s\in S$. \end{df}
Monoids are perhaps some of the most natural objects in mathematics. Take, for example, the natural numbers. One can add two natural numbers, and the number $0$ doesn't change what it is added to.
\begin{ex}
A slightly more complex example of a monoid is the set of all homeomorphism classes of compact oriented surfaces. The monoidal operation is the connected sum: given two manifolds $M$ and $N$, pick two small open sets homeomorphic to the open disc $D^2$ on $M$ and $N$ and form $M \# N = (M\setminus D^2)\cup_{\partial D^2}(N\setminus D^2)$. Intuitively, the picture is to cut two small discs out of $M$ and $N$ and glue the two surfaces together at their new boundary. Since $S^2\setminus D^2 = cl(D^2)$, the closure of $D^2$, taking any manifold $M$ and forming its connected sum with $S^2$ will produce a manifold homeomorphic to $M$. This example is in fact no different from the natural numbers $\mathbb{N}$ since if $M_g$ is the unique homeomorphism class of a closed genus $g$ oriented surface, $M_g \# M_h = M_{g+h}$. Both this monoid and the natural numbers $\mathbb{N}$ happen to be commutative but this is not always the case. \end{ex}
In a similar vein, certain categories will admit constructions such as the direct sum or tensor product which will ``feel'' monoidal. This idea is formalized in the definition for monoidal categories:
\begin{df} A ``monoidal category'' is a category $\mathcal{C}$ equipped with a bi-functor (functorial in both variables) $\otimes:\mathcal{C}\times \mathcal{C}\rightarrow \mathcal{C}$ such that $\otimes$ is associative up to natural isomorphisms that satisfy ``coherence conditions'' (for sake of brevity, it means that every commutative diagram that one would want to commute does) and there is a distinguished object $1\in\mathcal{C}$ such that $1\otimes X = X\otimes 1 = X$ for all $X\in\mathcal{C}$. The category is said to be ``symmetric monoidal'' if there is further a natural isomorphism $A\otimes B\rightarrow B\otimes A$ which squares to the identity map and is compatible with the coherence maps of $\mathcal{C}$. \end{df}
\begin{ex} The category of modules over a commutative ring $R$ is a monoidal category with either $\oplus$ or $\otimes$ as its monoidal operation. In the first case, the unit is the trivial module $0$, in the second case, the unit is the free $R$-module $R$. Both these operations turn $R$-Mod into a symmetric monoidal category. If the ring $R$ is not commutative, however, then the category of $R$-$R$ bimodules is still a monoidal category under $\otimes$, but this operation is no longer symmetric.\end{ex}
The symmetric monoidal category that this thesis will primarily be concerned with is the cobordism category $\text{Cob}(n)$.
\begin{df} Two $n$ dimensional manifolds $M$ and $N$ are said to be ``cobordant'' if there is an $n+1$ dimensional manifold $B$ with boundary $\partial B=M\coprod N$. The manifold $B$ is called a ``bordism.'' The ``cobordism category'' $\text{Cob}(n)$ is the category that has diffeomorphism classes of $n$ dimensional compact oriented manifolds as objects, and bordisms between manifolds are morphisms. More explicitly, $B:M\rightarrow N$ is a morphism if $\partial B = \bar{M}\coprod N$ where $\bar{M}$ is $M$ with the opposite orientation. Composition of bordisms $B:L\rightarrow M$ and $B':M\rightarrow N$ is given by gluing $B$ to $B'$ along $M$. The cobordism category can be made into a symmetric monoidal category with the disjoint union as its monoidal operation and the empty set (considered as an $n$ manifold) is the unit. \end{df}
Notice that the category $\text{Cob}(n)$ exhibits an unusually high amount of duality: given an $n+1$ dimensional manifold $B$, any partition of its boundary pieces into two disjoint sets gives a new morphism. For example, consider the circle $S^1\in\text{Cob}(1)$. The cylinder $S^1\times I$ can be thought of as three different bordisms: $B:S^1\rightarrow S^1$, $B:\emptyset \rightarrow \bar{S}^1\coprod S^1$, and $B:S^1\coprod \bar{S}^1\rightarrow \emptyset$. It turns out this flexibility in $\text{Cob}(n)$ will be instrumental in studying the cobordism category.
\begin{df} A ``strict monoidal functor'' between two monoidal categories $(\mathcal{C},\otimes_{\mathcal{C}})$ and $(\mathcal{D},\otimes_{\mathcal{D}})$ is a functor $F:\mathcal{C}\rightarrow \mathcal{D}$ such that there are equalities $F(X \otimes_\mathcal{C} Y) = F(X)\otimes_{\mathcal{D}} F(Y)$. The functor $F$ is said to be symmetric if $\mathcal{C}$ and $\mathcal{D}$ are symmetric monoidal categories. \end{df}
\begin{ex} Consider the symmetric monoidal category of CW-complexes with the monoidal operation being the wedge product. Then homology $H^n$ is a symmetric monoidal functor to the category of abelian groups with direct sum as the monoidal operation. More generally, if $R$ is a commutative ring, then $H^n(-;R)$ is a symmetric monoidal functor with values in $R$-modules.\end{ex}
\section{Low Dimensional Topological Quantum Field Theories}
\begin{df} An ``$n+1$ dimensional topological quantum field theory'' is a strict monoidal functor $Z:\mathrm{Cob}(n)\rightarrow \mathrm{Vect}(k)$ where $\mathrm{Vect}(k)$ is the category of vector spaces over a field $k$ (in practice we will take $k=\mathbb{C}$) \cite{lurie09}. \end{df}
More concretely, this means to each compact $n$-manifold without boundary $M$ we assign $Z(M)$, a $k$-vector space, and to each cobordism $F:M\rightarrow N$ we assign a linear map $Z(F):Z(M)\rightarrow Z(N)$. Recall that a morphism in the cobordism category is just an oriented manifold $F$ such that $\partial F = \bar{M}\coprod N$. In particular, given a $n+1$ manifold $F$ without boundary, we can regard $F$ as a morphism $F:\emptyset\rightarrow \emptyset$ which means applying $Z$ gives $Z(F):k\rightarrow k$. Since all maps of this form must be scalar multiplication, we can think of $Z(F)$ as giving an element in $k$. Indeed, one already sees that given an $n+1$ dimensional topological quantum field theory, one can produce invariants of oriented $n+1$ manifolds $M$ by evaluating $Z(M)$.

Already mentioned earlier was how maps in $\mathrm{Cob}(n)$ can be interpreted in quite a few different ways. In particular, suppose one is given an $n+1$ manifold $B$ such that $\partial B = \bar{M}\coprod N$. Then one can think of $B$ as any of the following:
\begin{enumerate}
	\item{$B:M\rightarrow N$}
	\item{$B:\bar{M}\coprod N\rightarrow \emptyset$}
	\item{$B:\emptyset\rightarrow N\coprod \bar{M}$}
\end{enumerate}
\begin{figure}
\centering
\includegraphics[scale=1]{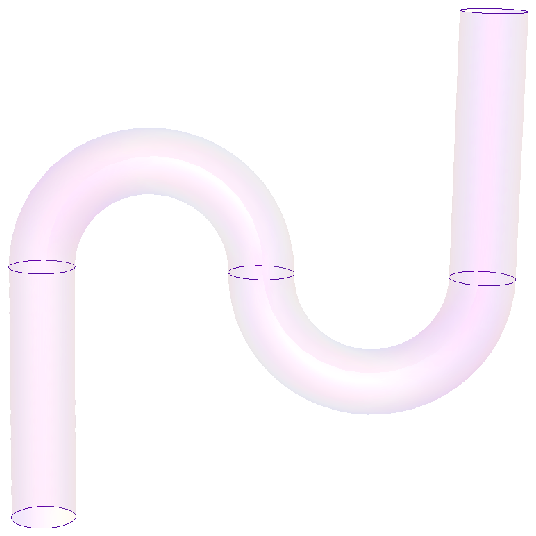}
	\caption{\label{fig:s}The composition $ev_M\coprod 1_M\circ 1_M\coprod coev_M$.}
\end{figure}
Letting $B=M\times [0,1]$ gives three maps: interpreted in the case of (1), $B:M\rightarrow M$ is the identity map; interpreted as (2), $B:\bar{M}\coprod M\rightarrow \emptyset$ gives evaluation; interpreted as (3), $B:\emptyset\rightarrow M\coprod \bar{M}$ gives the coevaluation map. This flexibility gives us a powerful tool in studying properties of the invariants that arise from a topological quantum field theory.
\begin{prop}
\label{DualOrientationsAreDualSpaces} Let $M$ be an oriented $n$-manifold and $Z$ an $n+1$ dimensional topological field theory. Then $Z(\bar{M}) = Z(M)^\vee$ where $\vee$ indicates taking the dual vector space. \end{prop}
\bpf
Consider the composition:
\[ M\stackrel{1_M\coprod coev_M}{\longrightarrow} M\coprod \bar{M}\coprod M\stackrel{ev_M\coprod 1_M}{\longrightarrow} M \]
If one were to draw a picture corresponding to this composition of bordisms, it would look similar to figure \ref{fig:s}. One can see that the composite $n+1$ manifold looks much like a S-shaped cylinder on $M$. Therefore, one can simply stretch this cylinder out and see it is clearly diffeomorphic to $M\times [0,1]$ hence the composition is the identity. This means that the induced linear map $Z((ev_M\coprod 1_M)\circ (1_M\coprod coev_M))$ must be the identity map on $Z(M)$. This condition is in fact very strong and generalizations of this condition will return later to control the behavior of topological quantum field theories (see example \ref{VectorSpaceDuals}). Now, the claim is that $ev_M:Z(M)\otimes Z(\bar{M})\rightarrow k$ induces a perfect pairing between $Z(M)$ and $Z(\bar{M})$. Suppose that $ev_M$ is degenerate; this means that there is a $w\neq 0$ such that for $v(x,w)=0$ for all $x\in Z(M)$. Looking at the image of $w$ gives:
\[ w\neq 0\mapsto 1\otimes w\mapsto (w'\otimes x)\otimes w = w'\otimes(x\otimes w)\mapsto w'\otimes v(x,w) = 0 \]
Therefore, this map cannot possibly compose to the identity and one gets that $v(x,w)=0$ for all $x$ implies $w=0$. Therefore, $v$ is non-degenerate and induces a perfect pairing which means $Z(M)$ is finite dimensional and $Z(\bar{M}) = Z(M)^\vee$. \epf
\subsection{Classificaion of Topological Quantum Field Theories in $1$ dimension}
Suppose $Z:\mathrm{Cob}(0)\rightarrow \mathrm{Vect}(k)$ is a topological quantum field theory. These theories are quite easy to define since the category $\mathrm{Cob}(0)$ only consists of two objects (and their disjoint unions): the positively and negatively oriented point denoted $+$ and $-$ respectively. By theorem \ref{DualOrientationsAreDualSpaces}, we see that $Z(+)=Z(-)^\vee$ so these vector spaces must be finite dimensional. Since in general any oriented $0$ manifold is just a set of positively oriented points and a set of negatively oriented points, working out the morphisms from a single point to a single point will determine the entire topological quantum field theory (notice there is no bordism from a single point to two points; such a bordism would look like a Y shape which is not a $1$-manifold at junction point). Any bordism $B$ between two points $+$ and $-$ is essentially just a cylinder and therefore is the identity map on $Z(+)$. Equivalently, it also gives the coevaluation map $k\rightarrow Z(-)\otimes Z(+)$ and evaluation $Z(+)\otimes Z(-)\rightarrow k$. Importantly, all of these maps are determined by the choice of $Z(+)$. If $B$ is a $1$-manifold without boundary, then $B$ must be diffeomorphic to disjoint unions of $S^1$. Therefore, it suffices to calculate the value of $Z(S^1)$. To do this, we break up $S^1$ into two half-circles connecting $-$ to $+$. Then $Z(S^1)$ is just the composition $ev\circ coev$; since $Z(+)$ is a finite dimensional vector space, the evaluation and coevaluation maps are easy to explicitly describe.
\begin{lem} If $V$ is a finite dimensional vector space, then there is a canonical isomorphism $V^\vee\otimes V\rightarrow \hom(V,V)$. \end{lem}
\bpf We define the map $V^\vee\otimes V\rightarrow \hom(V,V)$ to be the following: $(v_i^\vee, v_j)\mapsto (v_k\mapsto (v_i^\vee v_k)v_j)$. Since $V$ is finite dimensional, $\dim(\hom(V,V))=\dim(V^\vee\otimes V)=\dim(V)^2$. Now, pick a basis $\{v_i\}$ of $V$; we automatically are given a dual basis $\{v_i^\vee\}$ of $V^\vee$. Then $(v_i^\vee v_k)v_j$ is equal to $v_j$ whenever $i=k$ and $0$ everywhere else. These maps clearly form a basis for $\hom(V,V)$ therefore the map has maximal rank and is an isomorphism. \epf
Now, the coevaluation map can be described as $coev(r)=r1_V$ for $r\in k$ and $1_V$ the identity map in $\hom(V,V)$. Dually, the evaluation map is $ev(T) = tr(T)$ for some $T\in \hom(V,V)$. Therefore, the image of $1$ under the composition is $tr(1_V)=\dim(V)$. Therefore, the data of a topological quantum field theory reduces to choosing a finite dimensional vector space $V$ and $Z(S^1)=\dim(V)$.
\subsection{Classification of Topological Quantum Field Theories in $2$ dimensions}
\begin{figure}
\centering
\includegraphics[scale=.6]{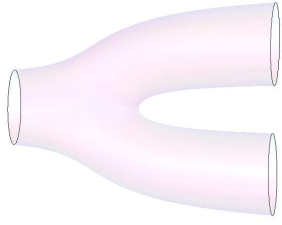}
	\caption{\label{fig:pp}The ``pair of pants.''}
\end{figure}

Topological quantum field theories in $2$ dimensions are only slightly more complicated than their $1$ dimensional counterparts. This is largely because the category $\mathrm{Cob}(1)$ remains an easily describable category: the objects are compact oriented $1$ dimensional manifolds which are just disjoint unions of circles. Therefore, the object assigned to $Z(\coprod_{i\in I}S^1)$ is controlled by the finite dimensional vector space $V=Z(S^1)$. A bordism $B:\coprod_{i\in I}S^1 \rightarrow \coprod_{j\in J}S^1$ is just a ``generalized pair of pants'' - a pair of pants with one waist hole for each $i\in I$ and one leg hole for each $j\in J$. Since both $I$ and $J$ are finite sets, as required by compactness, one can break such a bordism down into a composition of bordisms going from one circle to two circles of vice versa (figure \ref{fig:pp}). Therefore, one only has to consider what objects to assign to $Z(S^1)$ and what the bordism $Z(S^1\coprod S^1)\rightarrow Z(S^1)$ looks like. In particular, since $Z(S^1\coprod S^1)= V\otimes V$, one has a multiplication map $Z(S^1 \coprod S^1)=V\otimes V\stackrel{\mu}{\longrightarrow} V=Z(S^1)$.
\begin{figure}
\centering
\includegraphics[scale=.6]{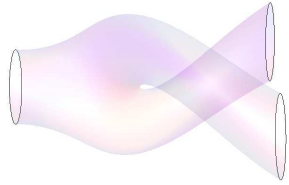}
		\caption{\label{fig:twistedpp}The commutativity of multiplication.}
\end{figure}
By drawing out the relevant ``pants'' diagrams, one can see that $\mu$ must be commutative and associative (figures \ref{fig:twistedpp} and \ref{fig:assocpp}). Now, consider the closed disc $D$ as an oriented $2$-manifold with boundary. One can view $D$ in two ways: $D:\emptyset\rightarrow S^1$ and $D:S^1\rightarrow \emptyset$. At first this may seem disconcerting because we do not seem to be tracking the orientation of $S^1$. This is in fact okay because there is a diffeomorphism $S^1\rightarrow \bar{S^1}$ (for example, taking $(x,y)\in\mathbb{R}^2$ to $(x,-y)$). The first interpretation of $D$ gives the ``unit'' map $u:k\rightarrow V$ and the second interpretation of $D$ gives the counit or ``trace'' map $tr:V\rightarrow k$. Indeed, $u(1)$ is the unit of the multiplication map $\mu$ since the corresponding composition of bordisms is diffeomorphic to the identity map on $S^1$. Now, if one multiplies then takes trace, the corresponding diagram is a half torus. This can be seen as the evaluation map $ev_{S^1}$ and gives a bilinear form $V\otimes V\stackrel{\mu}{\rightarrow}V\stackrel{tr}{\rightarrow}k$. Moreover, we checked earlier that this bilinear form gives a non-degenerate pairing. These commutative algebras with non-degenerate bilinear forms are familiar objects in mathematics:
\begin{figure}
\centering
\includegraphics[scale=1]{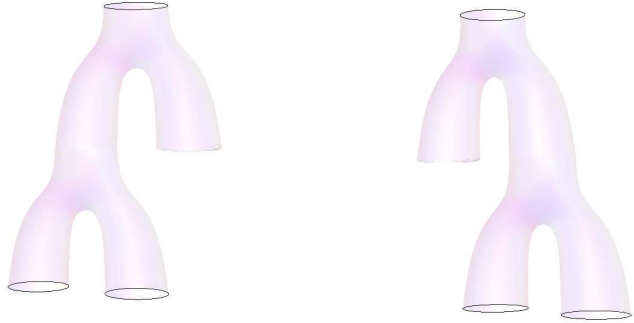}
		\caption{\label{fig:assocpp}The associativity of multiplication.}
\end{figure}
\begin{df} A ``Frobenius algebra'' is an associative unital algebra $A$ over $k$ and a non-degenerate bilinear form $T:A\times A\rightarrow k$. \end{df}
Therefore one sees that given a $2$ dimensional topological quantum field theory, one can construct a finite dimensional commutative frobenius algebra $A$. Conversely, suppose we are given such a Frobenius algebra: we may extract a functor $Z:\mathrm{Cob}(1)\rightarrow \mathrm{Vect}(k)$ as follows:
\begin{enumerate}
	\item{The object $Z(S^1)$ is just the underlying finite dimensional vector space of $A$}
	\item{The generating pair of pants diagram $B:S^1\coprod S^1\rightarrow S^1$ is given by the multiplication map $\mu:A\otimes A\rightarrow A$}
	\item{The unit map is just $u(1)=1_A$ and the counit $tr(v)$ is just given by $T(x,y)$ for any $xy=v$}
\end{enumerate}
As a result, one has the following equivalence
\begin{thm}\emph{(Classification of $2$ dimensional topological quantum field theories)} There is a bijection
\[ \mathrm{Fun}^\otimes(\mathrm{Cob}(1),\mathrm{Vect}(k))\stackrel{\sim}{\longrightarrow} \mathrm{CFrob}(k) \]
between equivalence classes of symmetric monoidal functors from $\mathrm{Cob}(1)$ to $\mathrm{Vect}(k)$ and isomorphism classes of commutative finite dimensional Frobenius algebras \cite{lurie09}\cite{teleman12}.\end{thm}
Now one might wonder if there is some $2$ manifold $\Sigma$ that calculates $\dim(A)$ just as we had in $1$ dimension. In $1$ dimension, the construction was to take coevaluation and then evaluation so we try the same thing again: the resulting manifold $\Sigma$ is the torus. The coevaluation map takes $1\in k\mapsto 1_A\in \mathrm{End}(A)$. Therefore, $tr(1_A)=\dim(A)$. If we look at what element we have immediately prior to taking trace, we have a way of selecting an element $\omega_A = \mu(coev(1))$. In \cite{abrams00}, Abrams shows this element to be of the form:
\[ \omega_A = \sum_i e_i e_i^\vee \]
where $e_i^\vee$ is the dual basis to $e_i$. This distinguished ``characteristic'' has many interesting properties:
\begin{thm} For a commutative Frobenius algebra $A$, the following hold:
\begin{enumerate}
	\item{$A$ is an Artinian ring (it satisfies the descending chain condition on ideals).}
	\item{If  $B$ is another commutative Frobenius algebra, then $\omega_{A\oplus B} =
	\omega_A\oplus \omega_B$.}
	\item{The minimal ideal of $A$ that intersects all other ideals of $A$ is $(\omega_A)$.}
	\item{$\omega_A$ is a unit if and only if $A$ is semisimple.}
	\item{A closed genus $g$ surface $M_g$ computes $Z(M_g)=tr(\omega^g)$.}
\end{enumerate}
(See \cite{teleman12} and \cite{abrams00}).
\end{thm}
\subsection{A Look Towards Extended Topological Quantum Field Theories}
While topological quantum field theories are described nicely as symmetric monoidal functors $Z\in \mathrm{Fun}^\otimes(\mathrm{Cob}(n),\mathrm{Vect}(k))$, these functors require increasingly more information to define for larger $n$. In the case of $n=0,1$, the manifolds were easy to describe and therefore these theories were controlled by a relatively small amount of data. One potential remedy for this increasing complexity is if one were able to build higher dimensional manifolds out of lower dimensional ones. For example, a cylinder $S^1\times I$ might be viewed as the identity bordism $S^1\rightarrow S^1$ and then the boundary $S^1$ might be viewed as an endomorphism of a point. Heuristically, one can think of this as enhancing the categorical structure of $\mathrm{Cob}(n)$: the objects are points, morphisms are bordisms between points, and ``2-morphisms'' are bordisms between bordisms of points and so on. If we then require the functor $Z$ to assign invariants to manifolds and morphisms of higher codimension, one can hope that the functor $Z$ might satisfy more rigid constraints allowing them to be easier to construct or classify. In fact, we've already seen one example of a topological field theory which was truly easy to describe: the $1$ dimensional theories. In this case, one completely formally obtains duality conditions on the point - namely that the evaluation and coevaluation satisfy a pairing. Using a diagram in $\mathrm{Cob}(0)$ we concluded that the theory was determined by evaluating $Z(+)$ and that the resulting vector space had to be finite dimensional with $Z(-)=Z(+)^\vee$, a kind of dualizability condition. These ideas are formalized in the study of ``extended topological quantum field theories'' which roughly speaking, are theories that assign invariants to manifolds of dimension less than $n$. In the best case scenario, if the theory can be ``extended down to points,'' then one can guess that $Z(+)$ might control a large amount of the theory and that only a select special class of objects might be assigned to $Z(+)$. In fact, this best case scenario is true: if one has a fully extended topological quantum field theory, then it is completely determined by $Z(+)$ and one can write down formal categorical conditions that $Z(+)$ must satisfy. The main difficulty is that as we include more and more dimensions, these categorical duality conditions become increasingly complex. However, the extended $2$ dimensional case still remains a tractable problem and it is easy to identify exactly what kinds of objects correspond to extended $2$ dimensional theories. Before proceeding to classify the extended $2$ dimensional theories, one must make rigorous exactly what kinds of ``higher categories'' $\mathrm{Cob}(n)$ and $\mathrm{Vect}(k)$ should be replaced with.

\section{Higher Categories}
If one of the major points of interests of topological quantum field theory is its use in identifying interesting invariants of manifolds, one major flaw is the inability to calculate these invariants using smaller pieces of the manifold. Since an $n+1$ dimensional topological quantum field theory does not give invariants of manifolds of dimension less than $n-1$, in general, one is not able to calculate invariants for smaller objects and combinatorially glue them together. For example, if one triangulates an $n+1$ dimensional manifold $M$, there is no procedure to obtain $Z(M)$ from the simplicial data of points, edges, faces, etc.

Ideally, one can enrich the domain of our functor $Z$ so that $Z$ accepts manifolds of arbitrary dimension less than $n$ as an argument while producing invariants in a sensible way. Furthermore, it turns out that if such an ``extension'' of $Z$ exists, the parameters required to define $Z$ become highly constrained ultimately allowing for a classification theorem. As such, we will now devote some time into understanding how to enrich the source and target categories of $Z$. Throughout this chapter, we will mostly follow Lurie's exposition in \cite{lurie09}.
\subsection{Strict 2-Categories}
Before jumping into the theory of $n$-categories, it is instructive to first look at what one expects out of a theory of higher categories in the simplest case. Recall that one of our goals was to be able to have $Z$ produce invariants of manifolds of all dimensions. Therefore, whatever our source category is, it in some sense should have many ``tiers'' of structure with each tier mimetic of $\mathrm{Cob}(n)$ for some $n$. With this in mind, we can attempt our first definition of a $2$-category.
\begin{df} A ``strict $2$-category'' $\mathcal{C}$ consists of the data:
	\begin{enumerate}
		\item{A class of objects $Ob(\mathcal{C})$.}
		\item{For any two objects $x$ and $y$, a category $Maps_{\mathcal{C}}(x,y)$ (we will omit the subscript when the category is clear from context).}
		\item{For each object $x$, a distinguished object $1_x\in Maps(x,x)$ called the identity map.}
		\item{A bifunctor called ``composition'' $\circ: Maps(x,y)\times Maps(y,z)\rightarrow Maps(x,z)$ such that composing with $1_x$ on the left and right does not change the morphism}
		\item{The composition bifunctor is strictly associative}
	\end{enumerate}
\end{df}
\begin{df} A ``functor between $2$-categories'' $F:\mathcal{C}\rightarrow \mathcal{D}$ is a rule assigning each object of $\mathcal{C}$ an object of $\mathcal{D}$ given any $x,y\in \mathcal{C}$, a functor $F:Maps_\mathcal{C}(x,y)\rightarrow Maps_\mathcal{D}(Fx, Fy)$ such that composition and identities are respected.\end{df}
Now, recall that for any given object $x$ in an arbitrary category, one can always form the endomorphism monoid $\mathrm{End}(x)$. Elements in this monoid are just endomorphisms of $x$ and the multiplicative structure is simply composition of endomorphisms. In the case of higher categories, one would expect that there is a similar, ``categorified'' version of $\mathrm{End}(x)$:
\begin{prop} Let $\mathcal{C}$ be a $2$-category and $x\in\mathcal{C}$. Then the category $\Omega_x\mathcal{C}=Maps_{\mathcal{C}}(x,x)$ is a monoidal category.\end{prop}
The notation here is inspired from topology where given a distinguished basepoint $x$ in a topological space $X$, one can form the loop space $\Omega_x X$ of loops based at $x$. The analogy between topological spaces and higher categories will be elaborated on in section \ref{ShortAside}.
\bpf The proof of this proposition is essentially a review of the definition of a $2$-category. To begin, one already knows that $\Omega_x\mathcal{C}$ is a category. Now, the monoidal bifunctor $\circ:\Omega_x\mathcal{C}\times\Omega_x\mathcal{C}\rightarrow \Omega_x\mathcal{C}$ is just the composition of morphisms. Associativity follows from the strict associativity of $2$-categories and the unit of composition is just $1_x$. \epf
Conversely, given any monoidal category $\mathcal{D}$, one can ``deloop'' the category into a strict $2$-category:
\begin{prop} Let $\mathcal{D}$ be a monoidal category. Define the category $B\mathcal{D}$ as follows: there is one object of $B\mathcal{D}$ which we will call the basepoint $\ast$ and $Maps_{B\mathcal{D}}(\ast,\ast) = \mathcal{D}$. Then $B\mathcal{D}$ is a strict $2$-category called the ``one-object delooping of $\mathcal{D}$.''
\end{prop}
\bpf The proof for this proposition is in the same spirit as above; since there is only one object $\ast$, there is only one mapping category $Maps(\ast,\ast) = \mathcal{D}$ which is given to us by hypothesis. Therefore the only remaining thing to check is strict associativity. However, since $\ast$ is the only object in $B\mathcal{D}$, strict associativity reduces to associativity in $\mathcal{D}$ which is given to us since $\mathcal{D}$ was monoidal. Finally, the unit of the monoidal structure is just the unit of $\mathcal{D}$ now treated as the morphism $1_\ast$. \epf
Later on, these constructions will make an apperance again in order to generalize notions of ``dualizability.''

We now examine some basic examples of $2$-categories:
\begin{ex} Let $k$ be a field. The $2$-category $\text{Vect}_2(k)$ of $k$-linear categories has:
	\begin{enumerate}
		\item{Objects are cocomplete categories enriched over $\text{Vect}(k)$.}
		\item{$1$-morphisms are cocontinuous $k$-linear functors.}
		\item{$2$-morphisms between two functors $F$ and $G$ are natural transformations from $F$ to $G$.}
		\item{Composition of maps is just composition of functors.}
		\item{The unit of composition is the identity functor $1_\mathcal{C}$.}
	\end{enumerate}
\end{ex}
Notice that given any unital $k$-algebra $A$, one can form a $k$-linear category by looking at left or right modules over $A$. Since module categories are cocomplete and complete, we will later only look at objects in $\text{Vect}_2(k)$ that arise as module categories as cocontinuous functors between module categories are easy to describe. This $2$-category is referred to as $\text{Alg}_2(k)$ (see section \ref{dualizabilitycriteria}).
\begin{ex} When attempting to define a higher cobordism category in a similar way, one might run into problems. For example, a naive definition for $\text{Cob}_2(1)$, the $2$-category of bordisms of $1$-dimensional manifolds might be
	\begin{enumerate}
		\item{Objects are oriented points.}
		\item{Given two points $p$ and $q$, $Maps(p,q)$ is the category of $1$-manifolds with boundary $p\coprod q$ and bordisms between them; in other words, given two bordsisms $M_1:p\rightarrow q$ and $M_2:p\rightarrow q$, a $2$-bordism $B:M_1\rightarrow M_2$ is a manifold with corners such that $\partial B=M_1\coprod_{p,q}M_2$ is a $1$-manifold with corners at $p$ and $q$.}
		\item{Composition of bordisms is given by gluing bordsism along the boundary.}
		\item{The units of composition are cylinders.}
	\end{enumerate}
	While this definition sounds appropriate at first, there is a problem with the way composition is defined in this category: in particular, consider the task of gluing to manifolds $M_1$ and $M_2$ along a boundary manifold $X$; there is no ``single'' manifold that ``is'' the gluing of $M_1$ and $M_2$ along $X$. Instead, the gluing is only defined up to unique isomorphism by the following universal property:
	\begin{df} Let $M_1$ and $M_2$ be manifolds each with $X$ as a boundary component. If $M$ is a third manifold equipped with maps $M_1\rightarrow M$ and $M_2\rightarrow M$ such that there is a commutative diagram
	\[\xymatrix{
	X\ar[r]\ar[d] &M_1\ar[d] \\
	M_2\ar[r] &M
	}\]
	then $M$ is said to be the ``gluing of $M_1$ and $M_2$ along $X$'' if the diagram is a pushout square; in other words, for any other manifold $N$ equipped with maps $M_1\rightarrow N$ and $M_2\rightarrow N$ which agree on $X$, there is a unique map $M\rightarrow N$ making the diagram commute
	\[\xymatrix{
	X\ar[r]\ar[d] &M_1\ar[d]\ar[rdd] \\
	M_2\ar[r]\ar[rrd] &M\ar[rd] \\
	&& N}\]
	Any manifold $M$ satisfying this property is well defined up to unique isomorphism.	\end{df}
Therefore the problem with composition of morphisms is that the notion of ``strictly associative'' must be loosened to some notion of ``associative up to isomorphism.'' This introduces a new problem: suppose one wishes to compose two $2$-morphisms in $\text{Cob}_2(1)$. If we want to make sense of having this composition ``up to isomorphism,'' then we would need to have maps between $2$-morphisms as well; i.e. we would need to introduce $3$-morphisms. Then this same problem might arise at when composing $3$-morphisms requiring us to add higher morphisms ad infinitum. Therefore, in order to formally solve this problem, one should get the sense that there are morphisms of all degrees, but all of them past a certain finite degree only keep track of isomorphisms. This leads us to our brief discussion on $(\infty,n)$-categories.
\end{ex}
\subsection{A Short Aside on $(\infty,n)$-Categories}
\label{ShortAside}
A full exposition on the theory of $(\infty,n)$-categories is unfortunately outside the scope of this thesis. This section only aims to sketch the motivating principles of the theory of higher categories without doing much to formalize the intuition. If interested in a rigorous treatment of the subject, the reader is again referred to \cite{lurie09}. The goal of this section is to provide an intuitive feel for $(\infty,n)$-categories and to explain how these objects might help solve the ``associativity up to coherent isomorphisms'' problem in the previous section. To begin, we examine one of the most natural places where higher categories arise:
\begin{ex} Let $X$ be a suitably good topological space. For each natural number $n\geq 0$, define a category called the fundamental $n$-groupoid of $X$ written $\Pi_{\leq n} X$:
	\begin{enumerate}
		\item{Objects are the points of $X$.}
		\item{$1$-morphisms are paths in $X$.}
		\item{$2$-morphsims are homotopies of $1$-morphisms.}
		\item{$n$-morphisms are homotopies of $n-1$-morphisms.}
	\end{enumerate}
	Notice that $\pi_0 X$ is just the isomorphism classes of $\Pi_{\leq 1}X$ and that $\pi_1(X,x) = \mathrm{Aut}_{\Pi_{\leq1}X}(x))$. This construction is kind of a primordial example of an $n$-category; however, there is a special feature of this $n$-category: since paths are invertible (by running $t\in [0,1]$ backwards) and homotopies are equivalence relations, every morphism in this $n$-category is invertible. Categories in which all morphisms are invertible are called ``groupoids.''
\end{ex}
Now, one would like some sort of notion for when $\Pi_{\leq n}X$ captures all of the homotopy-theoretic data for $X$:
\begin{df} A topological space $Y$ is called a ``homotopy $n$-type'' if $\pi_k(Y) = 0$ for all $k > n$. Given any space $X$, one can construct a homotopy $n$-type $Y$ with the same homotopy groups as $X$ up to dimension $n$; furthermore, the (weak) homotopy type of $Y$ is uniquely determined in this construction. \end{df}
As one would expect, if $X\rightarrow Y$ is a weak equivalence of homotopy $n$-types, then there is an induced equivalence of categories $\Pi_{\leq n}X\rightarrow \Pi_{\leq n}Y$. In general, however, homotopy groups exist in all dimensions; it is actually rather rare to encounter spaces with only finitely many non-zero homotopy groups. With this in mind, it would be appropriate to think of an arbitrary good topological space $X$ (for example, a simplicial set) as a homotopy $\infty$-type. Then, if $f:X\rightarrow Y$ is a map of topological spaces that induces equivalences $\Pi_{\leq n} X \rightarrow \Pi_{\leq n} Y$ for all $n$, $f$ is a weak homotopy equivalence.

Therefore, one can imagine creating a category $\Pi_{\leq \infty} X$ which is an $\infty$-groupoid that contains all of the homotopy-theoretic data of $X$. In fact, there should be an equivalence
\[ \text{Topological Spaces up to Weak Equivalence} \rightarrow \infty\text{-groupoids}\]
The fact that $\infty$-groupoids are combinatorial structures that keep track of homotopy-theoretic data suggests the following definition for an $\infty$-groupod:
\begin{df} A ``$(\infty,0)$-category'' or ``$\infty$-groupoid'' is a Kan complex.\end{df}
Recall that if $X$ is a topological space, there is a functor $S$ called the ``total simplicial complex'' taking $X$ to a Kan complex $SX$. Since Kan complexes are the fibrant objects in the model category of simplicial sets, after passing to homotopy theory, these Kan complexes give (weak) homotopy types. Having defined $(\infty,0)$-categories, one might try to define $(\infty,n)$-categories inductively.
\begin{df} An ``$(\infty,n)$-category'' is a category enriched over $(\infty,n-1)$-categories. Heuristically, an $(\infty,n)$-category has $m$-morphisms for all natural numbers $m$ such that if $f$ is an $m$ morphism with $m > n$ then $f$ is invertible.\end{df}
This inductive definition actually does give the correct definition for $(\infty,n)$-categories but one must be careful with what ``enriched'' means. Again, trying to mimic strict associativity in ordinary category at all levels in an $\infty$-category can raise issues like in the case of the extended cobordism category. Essentially, the issue boils down to the following: if $X$ is a topological space such that precomposing with the Hopf fibration gives a nontrivial map $\pi_2X\rightarrow \pi_3X$, then the homotopy $3$-type of $X$ cannot be modelled as a strict category; if it were strict, it turns out the groupoid modelling it would correspond to an infinite loop space which should not be the case. See \cite{camarena} for example.

Because $(\infty, n)$-categories are able to keep track of morphisms of all degrees, these $\infty$-categories are perfect candidates for providing a way of tracking the higher isomorphisms in $\text{Cob}_2(1)$. Therefore, one should think of $\text{Cob}_2(1)$ as an $(\infty,2)$-category. Having defined higher categories, we can state the definition for a fully extended topological quantum field theory.
\begin{df}
Let $\mathcal{C}$ be a symmetric monoidal $(\infty,n+1)$-category. A (fully) extended $n+1$-dimensional topological quantum field theory is a symmetric monoidal functor
\[ Z: \mathrm{Cob}_{n+1}(n)\rightarrow \mathcal{C}\]
Here, the category on the left is the higher analogue of the category $\text{Cob}_2(1)$ which has points as objects, bordisms of points as $1$-morphisms, bordisms of bordisms as $2$-morphisms up to $n$-morphisms. Then the remaining infinite tiers of maps are used to track coherence data of structural isomorphisms. Since the rest of this thesis mainly concerns itself with the case $n=2$ for which we may do without a developed theory of $\infty$-categories, we will not say more on this.
\end{df}
\subsection{Dualizability in Higher Categories}

Recall that when classifying $1$ dimensional topological quantum field theories, one must assign a vector space $V=Z(+)$ to a positively oriented point $+$. Then the duality between the positively oriented point and the negatively oriented point produced a perfect pairing between $V$ and $V^\vee$ inducing an isomorphism $V\cong (V^\vee)^\vee$ which implies $\dim(V)$ is finite. In general, if we wish to produce an extended topological quantum field theory, there will be some sort of generalized duality condition which comes from the duality of the point objects in the higher cobordism category. The higher dimension the topological quantum field theory, the stronger the duality conditions. We now introduce the higher dualizability constraints that $Z(+)$ must satisfy for extended topological quantum field theories.
\begin{df} Let $\mathcal{C}$ be a $2$-category. If $f:X\rightarrow Y$ and $g:Y\rightarrow X$ are two $1$-morphisms, then a $2$-morphism $u:1_X\rightarrow g\circ f$ is said to be the unit of adjunction if there exists a $2$-morphism $v:f\circ g\rightarrow 1_Y$ such that the following two compositions are the identity $2$-morphisms:
	\begin{eqnarray}
		&f = f\circ 1_X \stackrel{1_f\times u}{\longrightarrow}f\circ g\circ f\stackrel{v\times 1_f}{\longrightarrow}1_Y\circ f = f \nonumber \\
		&g = 1_X\circ g \stackrel{u\times 1_g}{\longrightarrow}g\circ f\circ g\stackrel{1_g\times v}{\longrightarrow}g\circ 1_Y = g \nonumber
	\end{eqnarray}
	In this situation, one says that $f$ is the left adjoint of $g$ and $X$ is the dual of $Y$. If every $1$-morphism in $\mathcal{C}$ admits both a left and right adjoint, then $\mathcal{C}$ is said to ``have adjoints.''
\end{df}
Now suppose $\mathcal{C}$ is just any ordinary monoidal category then $\mathcal{C}$ is said to ``have duals'' if $B\mathcal{C}$ has adjoints.
\begin{ex}
\label{VectorSpaceDuals}
Consider the monoidal category $\text{Vect}(k)$ with the operation $\otimes_k$. We can deloop this category into $B\text{Vect}(k)$ and ask when a vector space (now considered as a $1$-morphism) has adjoints. In particular, if $V:\ast\rightarrow \ast$ is dualizable, one must provide $W:\ast\rightarrow \ast$ along with units of adjunction $v:V\otimes W\rightarrow k$ and $u:k\rightarrow W\otimes V$ such that the following $2$-morphisms compose to the identity:
	\begin{eqnarray}
		&V = V\otimes k \stackrel{1_V\times u}{\longrightarrow}V\otimes W\otimes V\stackrel{v\times 1_V}{\longrightarrow} k\otimes V = V \nonumber \\
		&W = k\otimes W\stackrel{u\times 1_W}{\longrightarrow}W\otimes V\otimes W\stackrel{1_W\times v}{\longrightarrow}W\otimes k = W \nonumber
	\end{eqnarray}
	At this point, the example should appear familiar as it is identical to the proof of proposition \ref{DualOrientationsAreDualSpaces}. As before, $v:V\otimes W\rightarrow k$ induces a perfect pairing between $V$ and $W$. From this, one sees that adjunction generalizes the duality exhibited in the case of non-extended topological quantum field theory (in other words, this is the correct notion of dualizability one layer down in the theory).
\end{ex}
\begin{df} An object $X$ in a symmetric monoidal $n$-category $\mathcal{C}$ is said to be dualizable if $X$ has adjoints considered as a $1$-morphism in $B\text{ho}_1\mathcal{C}$ where $\text{ho}_1\mathcal{C}$ is the ordinary category whose objects are the objects of $\mathcal{C}$ with morphisms isomorphism classes of $1$-morphisms in $\mathcal{C}$. The category $\text{ho}_1\mathcal{C}$ is called the ``homotopy category'' of $\mathcal{C}$. \end{df}
While we have not rigorously defined exactly what ``symmetric monoidal $n$-category'' means, in practice this will not be much of an issue since the primary case considered is the category $\text{Alg}_2(k)$ which has a intuitively clear tensor product: given two $k$-algebras $A$ and $B$ one can form the tensor product $A\otimes_k B$. Now, suppose one had a bimodule ${_B}M_A$. We need a way to produce a $B\otimes C$-$A\otimes C$ bimodule: to do this consider the module $M\otimes_k C$. Define the following actions:
\begin{eqnarray}
&& (b\otimes c)(m\otimes c') = (bm\otimes cc') \nonumber\\
&& (m\otimes c')(a\otimes c) = (ma\otimes c'c)\nonumber
\end{eqnarray}
This turns $M\otimes_k C$ into the required morphism $A\otimes C\rightarrow B\otimes C$. This is the approach taken by, for example, Schommer-Pries in \cite{schommerpries11}.

In general, to define the adjunction in higher $n$-categories (for $n > 2$), one has to do a little more work. Since this thesis will not use anything higher than a $2$-category, we will not go into detail, but the general idea is to form a $2$-category which only keeps track of $1$-morphisms and isomorphism classes of $2$-morphisms and then study adjunction in this $2$-category. Then one iteratively applies this construction to all the mapping categories until all the morphism levels of the category are exhausted
\begin{df} Let $\mathcal{C}$ be a symmetric monoidal $2$-category. Then there is a subcategory called the subcategory of ``fully dualizable objects'' whose objects consist only of objects in $\mathcal{C}$ that are dualizable and the only $1$-morphisms between objects are invertible morphisms that have both left and right adjoints. Objects in this category are said to be ``fully dualizable.''\end{df}
\begin{thm}
\label{Lurie2FullyDualizable}Let $\mathcal{C}$ be as above, then an object $X$ is fully dualizable if and only if the object $X$ admits a dual $X^\vee$ and the evaluation morphism $ev:X\otimes X^\vee\rightarrow 1$ has both left and right adjoints.\end{thm}
\bpf The proof for this theorem can be found in \cite{lurie09} section 4.2. The important point which we will take away is the factorization of the right and left adjoints of evaluation as
\begin{eqnarray}
&& ev^R_X=(S\otimes 1_{X^\vee})\circ coev_X \nonumber \\
&& ev^L_X=(T\otimes 1_{X^\vee})\circ coev_X\nonumber
\end{eqnarray}
It turns out that $S$ and $T$ are inverse isomorphisms are adjoints to one another; in fact, $S$ can be described as the Serre automorphism (constructed in section \ref{SerreAutomorphism}). This factorization will become useful as it will allow for the explicit construction of the Serre automorphism and will later give a ``trace map'' on fully dualizable objects.\epf

This notion of fully dualizable will be of utmost importance when studying extended topological quantum field theories; as it turns out, the condition that an object be fully dualizable is highly stringent. In the next chapter, we will formulate the ``Cobordism Hypothesis'' which will make apparent the importance of understanding fully dualizable objects and we will also calculate what the fully dualizable objects are in the category we earlier described as $\text{Alg}_2(k)$. \section{The Cobordism Hypothesis and the Classification of Extended 2D Topological Quantum Field Theories}
\subsection{The Cobordism Hypothesis}
The Cobordism Hypothesis is the main tool with which to study fully extended topological quantum field theories.
\begin{thm}
\emph{The Cobordism Hypothesis}
\label{CobordismHypothesis}
Let $\mathcal{C}$ be a symmetric monoidal $(\infty,n+1)$-category and $\mathrm{Cob}^{fr}(n)$ the $(\infty,n+1)$-category of framed bordisms. Then there is an $\infty$-groupoid $\mathcal{C}^{fd}$ obtained by taking only the fully dualizable elements of $\mathcal{C}$ and only considering the invertible morphisms with both left and right adjoints. The evaluation $Z\mapsto Z(\ast)$ of $Z$ at a point $\ast$ determines an equivalence of categories
\[ \mathrm{Fun}^\otimes(\mathrm{Cob}^{fr}(n), \mathcal{C}) \stackrel{\sim}{\longrightarrow} \mathcal{C}^{fd}\]
\end{thm}
\bpf The proof for the Cobordism Hypothesis is far out of the scope of this thesis. Roughly speaking, the idea is to show that the category $\text{Cob}^{fr}(n)$ is ``freely generated by the point.'' As in the case with other free objects, morphisms (in this case functors) out of the freely generated object are determined by what the it does on the generator, the point $\ast$. For more details, the interested reader is referred to \cite{lurie09} for the general case and \cite{schommerpries11} for the $2$ dimensional case.\epf
One observation is that in the past, we only required the manifolds to be oriented whereas now we want a framing. It turns out, however, that given a framed manifold $M$, one can specify additional data to reduce the condition down to only requiring $M$ to be equipped with an orientation.
\begin{df} Let $V$ be an $n$-dimensional real vector space. Recall that a ``framing'' $\xi$ of $V$ is an ordered basis $(v_1,\cdots v_n)$ of $V$. In other words, a frame is a choice of isomorphism $\xi:\mathbb{R}^n\rightarrow V$.
\end{df}
Clearly, if one has two framings $\xi$ and $\xi'$, one can write an invertible linear transformation $A\in GL_n(\mathbb{R})$ which carries the basis $\xi$ onto $\xi'$. Therefore, if $Fr(V)$ denotes the set of framings of $V$, there is a right action $Fr(V)\times GL_n(\mathbb{R})\rightarrow Fr(V)$ by $(\xi, A)\mapsto \xi A$.
\begin{df} An ``orientation'' of $V$ is an equivalence class of framings where two framings $\xi$ and $\xi'$ are considered equivalent if the transformation $A$ taking $\xi$ to $\xi'$ has positive determinant. \end{df}
Now suppose one were more restrictive: if $V$ has an inner product, then one can ask that the frame $\xi$ is orthonormal. In this case, instead of the full $GL_n(\mathbb{R})$ acting on $Fr(V)$, one only has the orthogonal group $O(n)$ acting by isometries. If we were to restrict our attention to these orthonormal frames, then two frames determine the same orientation if and only if they differ by an element of $SO(2)$.
\begin{df} Let $E$ be a vector bundle of rank $n$ over a manifold $M$. An ``orthonormal frame'' $\xi$ is a choice of frame $\xi_x$ of the fiber $E_x$ over a point $x$ such that these choices vary smoothly as $x$ varies on $M$. A framing of $M$ is an orthonormal frame of the tangent bundle $TM$. Since over each point, a framing $\xi_x$ is just a trivialization of the fiber $E_x$, a framing of $M$ is a trivialization of $TM$. \end{df}
If $Fr(M)$ is the collection of orthonormal framings of $M$, one sees that again there is a right action $Fr(M)\times O(n)\rightarrow Fr(M)$ where $n=\dim(M)$. This right action is just the action of $O(n)$ on the frame $\xi_x$ over each point $x\in M$.
\begin{df} If $G\subseteq O(n)$ is a subgroup, a ``$G$-structure'' on $M$ is an equivalence class of bundles $[\xi]$ where two bundles $\xi$ and $\xi'$ are equivalent if there is an element $g\in G$ such that $g\xi=\xi'$. \end{df}
\begin{ex} We quickly look at 3 subgroups of $O(n)$ that provide familiar $G$-structures:
\begin{enumerate}
	\item{If $G$ is the trivial group $\{1\}$, then a $G$-structure is an orthonormal frame $\xi$.}
	\item{If $G$ is the entire group $O(n)$, then since $O(n)$ acts transitively on the set of orthonormal frames, a $G$-structure is no additional data at all.}
	\item{If $G$ is $SO(n)$, the connected component of the identity, then an $SO(n)$-structure is a smooth choice of orientation on the fiber of each point which is an orientation on $M$.}
\end{enumerate}
\end{ex}
Therefore, one sees that the framed bordism case is actually more general than the oriented bordism case. There is in fact a way to pass from framed bordisms to oriented bordisms.
\begin{df} There is an action of $O(n)$ on the category $\text{Fun}^\otimes(\text{Cob}^{fr}(n), \mathcal{C})$ by $gZ(M)=Z(g(M))$ where $g\in O(n)$ is acting on the framing $\xi$ of $M$. \end{df}
Since by the cobordism hypothesis, there is a categorical equivalence \[\text{Fun}^\otimes(\text{Cob}^{fr}(n), \mathcal{C})\stackrel{\sim}{\rightarrow}\mathcal{C}^{fd}\]
one also gets an action of $O(n)$ on $\mathcal{C}^{fd}$. Furthermore, since $\mathcal{C}^{fd}$ is an $\infty$-groupoid, we can consider it as uniquely determining a homotopy type $X$ such that $\mathcal{C}^{fd} \simeq \Pi_{\leq \infty}X$.
\begin{thm} If $G\subseteq O(n)$ is a subgroup and $\mathrm{Cob}^G(n)$ is the category of bordisms with $G$-structure (for manifolds of dimension less than $n$, one looks at framings of the stablized tangent bundle $TM\oplus \mathbb{R}^{n-\dim(M)}$). Then there is a categorical equivalence
\[\mathrm{Fun}^\otimes(\mathrm{Cob}^{G}(n), \mathcal{C})\stackrel{\sim}{\rightarrow}(\mathcal{C}^{fd})^{hG} \]
Where the category on the right denotes the homotopy fixed points of $\mathcal{C}^{fd}$.\end{thm}
\bpf Again, the reader is referred to \cite{lurie09} for a proof. \epf
In particular, the category of oriented extended $2$ dimensional topological quantum field theories with values in $\mathcal{C}$ is equivalent to $(\mathcal{C}^{fd})^{hSO(2)}$. We now describe this action of $SO(2)$:
\begin{df}
\label{SerreAutomorphism}
	The unit interval $[0,1]$ can be interpreted as the identity bordism. In order for this to be a framed bordism in $\text{Cob}^{fr}(1)$, one must pick a framing of $T[0,1]\oplus \mathbb{R}$ such that the framings at $0$ and $1$ are the same (say equal to the standard framing of $\mathbb{R}^2$). This set actually has a composition operation by gluing intervals together: $[0,1] \cup [1,2] = [0,2]\simeq [0,1]$ with the framings at the end points identified. In fact, the set of all such framings up to homotopy can be identified by how many times the framing twists as $x$ ranges over the interval. In other words, the set of these framings forms a $\pi_1(SO(2))$-torsor where an integer $z$ induces an additional $z$ twists on the framing of $[0,1]$. Since $SO(2)$ acts on itself, one can choose the left invariant framing of $S^1$ to be the identity element to obtain an isomorphism with $\pi_1(SO(2))$. Now, suppose one were given an arbitrary $g\in SO(n)$. If there is a path $p$ connecting $1\in SO(2)$ to $g$, then one can always trivialize the framing along this path as long as it represents $0\in \pi_1(SO(2))$. However, if one looks at the action of the generator $\gamma\in \pi_1(SO(2))$, the framing cannot be trivialized and therefore $\gamma$ gives a nontrivial automorphism on the point $+$. This nontrivial automorphism is called the ``Serre automorphism'' and is denoted $S$.
\end{df}
\begin{thm}
\label{trivializationofS} If $\mathcal{C}^{fd}$ is the space of fully dualizable objects for a $2$-category $\mathcal{C}$, an invertible $2$-morphism identifying the Serre automorphism of a point $x$ with the identity map at $x$ is equivalent to giving $x$ the structure of a $SO(2)$ homotopy fixed point.\end{thm}
Before proving this, we briefly recall a fundamental result of obstruction theory:
\begin{thm}
Let $B$ be a simply connected CW complex, $X\rightarrow E\rightarrow B$ be a fibration and suppose $\sigma_n:B_n\rightarrow E$ is a section on the $n$-skeleton of $B$. If $D^{n+1}$ is an $n+1$-cell of $B$, then $\sigma_n|_{\partial B_{n+1}}:S^n\rightarrow E$. Projecting this down to $B$ takes the image to $D^{n+1}$ so $\sigma_n(S^n)$ really lies in the fiber $X$. Hence for each $n+1$ cell, $\sigma_n$ defines an element in $\pi_{n}X$. Since formal linear combinations of $n+1$ cells precisely give cellular cochains, this data specifies a cochain $C^{n+1}(B,\pi_n X)$. This cochain is actually a cocycle and $\sigma_n$ extends to $B_{n+1}$ precisely when this cocycle is $0\in H^{n+1}(B,\pi_n X)$ \cite{hatcherAT}.\end{thm}
We can now prove theorem \ref{trivializationofS} using this tool:
\bpf Giving an $SO(2)$ homotopy fixed point on $X=\mathcal{C}^{fd}$ is the same thing as a section of the Borel fibration $X\rightarrow E=(X\times ESO(2))/SO(2)\rightarrow BSO(2)=\mathbb{C}\mathbb{P}^\infty$. The space $\mathbb{C}\mathbb{P}^\infty$ has a nice filtration:
\[ \ast = \mathbb{C}\mathbb{P}^0\subset \mathbb{C}\mathbb{P}^1 \subset \mathbb{C}\mathbb{P}^2 \cdots
\] with $\mathbb{C}\mathbb{P}^\infty = \bigcup_n \mathbb{C}\mathbb{P}^n$. At each step of this filtration, one can ask how to glue down cells onto $\mathbb{C}\mathbb{P}^{n-1}$ to get $\mathbb{C}\mathbb{P}^n$. Since $H^*\mathbb{C}\mathbb{P}^\infty = \mathbb{Z}[x]$ the one generator in degree $2$, there is only one cell in each even dimension; this disc $D^{2n}$ then is glued to $\mathbb{C}\mathbb{P}^{n-1}$ by the Hopf fibration $S^{2n-1}\rightarrow \mathbb{C}\mathbb{P}^{n-1}$.

Now suppose we understand sections $\sigma_{n-1}$ on $\mathbb{C}\mathbb{P}^{n-1}$. Since $\pi_1\mathbb{C}\mathbb{P}^\infty =0$, obstruction theory tells us that we can extend the section over the unique cell in dimension $2n$ of $\mathbb{C}\mathbb{P}^n$ if and only if $\sigma_{n-1}: S^{2n-1}=\partial D^{2n}\rightarrow X$ is nullhomotopic as an element of $\pi_{2n-1}X$. In general, this problem is not easily approachable, but in this specific case the lifting problem is rather simplified. In dimension $0$, we can always lift the basepoint $\ast\in\mathbb{C}\mathbb{P}^\infty$ to some point $x\in X$. The $SO(2)$ action then draws out a loop $\sigma_0:S^1\rightarrow X$ and it is precisely this loop that must vanish in $\pi_1X$ in order to extend this section. Then, the next condition that must be satisfied is a nullhomotopy condition in $\pi_3 X$. But this comes for free: since $X$ is the space of fully dualizable objects in a $2$-category, there are no morphisms of degree higher than $2$ so $\pi_m X = 0$ for all $m\geq 2$ (in other words, since $\mathcal{C}^{fd}$ is a $2$-groupoid, it comes as $\Pi_{\leq 2}X$ for a homotopy $2$-type $X$). Therefore, as long as the section defines a nullhomotopic loop in $\pi_1 X$, we can immediately lift $\sigma_0$ to all of $\mathbb{C}\mathbb{P}^\infty$. Since this loop $\sigma_0:S^1\rightarrow X$ is precisely the loop drawn out by the action of $SO(2)$, we know that the loop corresponds to the Serre automorphism $S:x\rightarrow x$. Therefore, a nullhomotopy is precisely an invertible $2$-morphism trivializing $S\simeq 1_x$. \epf

\subsection{Criteria for Full Dualizability}
\label{dualizabilitycriteria}
In this section, we will derive a set of necessary and sufficient conditions for an object to be dualizable in $\text{Alg}_2(k)$.
\begin{df} Recall that the $2$-category of ``algebras bimodules and intertwiners'' $\text{Alg}_2(k)$ has unital $k$-algebras as objects. Given two algebras $R$ and $S$, and the mapping category of $1$-morphisms $Maps_{\text{Alg}_2(k)}(R,S)$ is the category of cocontinuous functors from $R$-Mod to $S$-Mod. The $2$-morphsims are natural transformations of functors.\end{df}
As mentioned earlier, cocontinuous functors between module categories are actually easy to describe: the Eilenberg-Watts theorem is precisely the tool which we will use to describe these functors:
\begin{thm}
\emph{Eilenberg-Watts Theorem}
\label{EilenbergWatts} Let $R$ and $S$ be arbitrary rings with unit. Then there is an equivalence of categories between the category of $S$-$R$ bimodules and the functor category of cocontinuous functors from left $R$-modules to left $S$-modules.
\[ S\text{-}\mathrm{Mod}\text{-}R \stackrel{\sim}{\longrightarrow} \mathrm{Fun}_k^{cocont}(R\text{-}\mathrm{Mod},S\text{-}\mathrm{Mod})\]
by taking the module $_SM_R$ to $_SM_R\otimes_R$ \cite{eilenberg60} \cite{watts60} \cite{ivanov12}.\end{thm}
Therefore, $\text{Alg}_2(k)$ has the equivalent formulation of having the mapping category $Maps_{\text{Alg}_2(k)}(R,S)$ be the category of $S$-$R$ bimodules with linear maps as $2$-morphisms. Thus, $\text{Alg}_2(k)$ is commonly called ``the $2$-category of algebras, bimodules and intertinwers:''
\begin{prop} Every algebra $A\in \text{Alg}_2(k)$ has a dual. \end{prop}
\bpf We first begin by describing the one object delooping of the homotopy category: $B\text{ho}_1\text{Alg}_2(k)$ has a unique object $\ast$ and $1$-morphisms $Map(\ast,\ast)$ algebras over $k$ with bimodules as $2$-morphisms. Any algebra gives a map $A:\ast\rightarrow \ast$ and so asking for a dual is providing an algebra $A^\vee:\ast\rightarrow \ast$ such that there is a unit of adjunction $u:k\rightarrow A^\vee\otimes A$ and $v:A\otimes A^\vee\rightarrow k$. The claim is that we can set $A^\vee$ to be the opposite algebra $A^o$ (the algebra has the same underlying elements as $A$ but with multiplication from right to left) and the unit and counit maps are $_{A\otimes A^o} A$ and $A_{A\otimes A^o}$ respectively. Then the adjunction condition reduces to showing that the following composition is identity:
\[ A = A\otimes k\stackrel{A\otimes A}{\longrightarrow}A\otimes A^o\otimes A\stackrel{A\otimes A}{\longrightarrow}k\otimes A=A \]
In other words, one must show that $(A\otimes A)\otimes_{A\otimes A^o\otimes A}(A\otimes A)\cong A$ as an $A$-$A$ bimodule. To do this, it is helpful to label all of the rings and the actions. We will do this with numerical subscripts:
\[ A_0 \stackrel{_0A_1\otimes _2A_3}{\longrightarrow}A_1\otimes A_2^o\otimes A_3\stackrel{_1A_2\otimes _3A_4}{\longrightarrow A_4} \]
Therefore, we wish to show that there is an isomorphism
\[ _0A_1\otimes _2A_3 \otimes_{A_1\otimes A_2\otimes A_3} {_1}A_2 \otimes _3A_4 \rightarrow {_0}A_4
\]
Once we write down the action of the ground ring, it becomes apparent that there are not many choices of a morphism that are multilinear: in particular, we see that if we require the rings with a right subscript to be adjacent to a ring with the same subscript on the left, we are forced to make the morphism $(w\otimes x\otimes y\otimes z)\mapsto wyxz$.

We first show that this map is indeed multilinear:
\begin{enumerate}
	\item{The map is clearly $A_0$ linear on the left and $A_4$ linear on the right}
	\item{Suppose $a_1\in A_1$. Then $a_1(w\otimes x\otimes y\otimes z) = (wa_1 \otimes x\otimes y\otimes z) = (w\otimes x\otimes a_1y\otimes z)\mapsto wa_1yxz$ which is multilinear}
	\item{A similarly trivial calculation shows that the map is multilinear with respect to the actions of $A_2$ and $A_3$.}
\end{enumerate}
Now, this map is clearly surjective as $(a\otimes1\otimes1\otimes1)\mapsto a$ for any $a\in A$. Now, suppose $(w\otimes x\otimes y\otimes z)\mapsto 0\in A$. Then this means $wyxz = 0$ in $A$. So take the each of the elements $w,x,y,z$ and interpreting them as actions of $A$ on $1$, we can pull all the terms to the first term in the tensor product and get $(w\otimes x\otimes z\otimes y) = (wyxz\otimes 1\otimes 1\otimes 1)$ and since the first term is $0$ in $A$, the tensor product is $0$ and we see that the map is injective and hence an isomorphism. Therefore, any algebra is dualizable with $A^o$ as its dual. \epf
\begin{prop}
\label{ConditionsForLeftAdjoints}
Let $A$ and $B$ be two unital $k$-algebras. Say $F$ is a functor from left $A$-modules to left $B$-modules of the form $F={_B}M_A\otimes_A (-)$ where ${_B}M_A$ is a $B$-$A$ bimodule. Then this map has a left adjoint if and only if $M$ is a finitely generated and projective as an $A$-module. \end{prop}
\bpf
If there is a left adjoint to this functor, we also know it must be cocontinuous and so the Eilenberg-Watts theorem tells us the functor is of the form ${_A}N_B\otimes_B (-)$. Therefore, by tensor-hom adjunction, ${_B}M_A\otimes_A (-) = \hom_A({_A}N_B,-)$ since they are both right adjoints of ${_A}N_B\otimes_B(-)$. We will show that ${_A}N_B$ has the conditions required and conclude that ${_B}M_A$ must also have the same conditions. First we argue that ${_A}N_B$ must be projective: suppose we have a surjection ${_A}X\twoheadrightarrow {_A}Y$ and a map ${_A}N\rightarrow {_A}Y$. We must show there is a lift:
\[\xymatrix{
					&{_A}X\ar@{->>}[d] \\
	{_A}N\ar[r]\ar[ru]^{\exists} &{_A}Y
}\]
This is equivalent to checking that the map ${_A}X\twoheadrightarrow {_A}Y$ induces a surjection $\hom_A({_A}N,{_A}X)\rightarrow \hom_A({_A}N,{_A}Y)$. However, we know that $\hom({_A}N, -) = M_A\otimes(-)$ so this is equivalent to showing that ${_B}M_A \otimes X\rightarrow {_B}M_A\otimes Y$ is a surjection. This is automatic since tensoring is right exact, therefore we know that $N$ is projective. It remains to show that $N$ is finitely generated; suppose the converse: write ${_B}N_A$ as the union of finitely generated submodules
\[ {_A}N_{B}=\bigcup_{\alpha} N_\alpha \]
where each $N_\alpha$ is an $A$-$B$ bimodule. Consider the identity element $1\in \hom_A({_A}N_B, \bigcup_\alpha N_\alpha)$. By adjunction, this is equivalent to giving an element $x\in {_B}M_A\otimes \bigcup_\alpha N_\alpha$. However, in this tensor product, we can surely write $x$ as the finite linear combination of pure tensors, so there is some finitely generated submodule ${_A}N'_B$ which this identity map factors through. Therefore ${_A}N_B$ must have also been finitely generated. Since ${_A}N_B$ is finitely generated and projective, it is the direct summand of a free module $A^m$. Now, ${_B}M_A\otimes(-) = \hom_A({_A}N_B, -)$ so putting the module $A$ in the variable position gives $_{B}M_A = \hom_A({_A}N_B,A)$, the dual module of ${_A}N_B$. Therefore, ${_B}M_A$ is the dual module to this summand of a free module and so it must be a summand of the dual free module. Consequently, ${_B}M_A$ is also finitely generated and projective as a right $A$-module.

Conversely, suppose ${_B}M_A$ is finitely generated and projective over $A$; then we can write it as a summand of a free module. Then picking ${_A}N_B$ to be the corresponding summand of the dual free module shows that there is a left adjoint to ${_B}M_A\otimes(-)$.
\epf
\begin{prop}
\label{ConditionsForRightAdjoints}
In the same setup as above if $B$ is Noetherian, the functor ${_B}M_A\otimes(-)$ has a cocontinuous right adjoint if and only if $M$ is finitely generated over $B$.
\end{prop}
\bpf
By tensor-hom adjunction, we automatically know that the functor ${_B}M_A\otimes(-)$ has a right adjoint of the form $\hom({_B}M_A, -)$. However, we do not know that this right adjoint is cocontinuous; in order to have cocontinuity, we need that $\hom({_B}M_A,-)$ preserves all small coproducts (and therefore all small colimits). If $M$ is finitely generated, then $\hom({_B}M_A,-)$ certainly preserves coproducts. If $B$ is Noetherian, then the converse is true \cite{head72}. In our case, the Noetherian condition on $B$ will be guaranteed since $B$ will be the ground field $k$.\epf

We now apply this theorem to the category $\text{Alg}_2(k)$ to obtain a list of criteria for when an algebra $A$ is fully dualizable. Recall that theorem \ref{Lurie2FullyDualizable} says that an algebra $A$ is dualizable if and only if a dual object exists and the evaluation map has adjoints on both sides. Since all algebras $A$ have $A^o$ as their dual, it suffices to check when the evaluation map has adjoints. At this point, we specialize to the case where $k$ is perfect, for example $k=\mathbb{C}$.
\begin{thm}
\label{FullyDualizableAlgebras} An algebra $A\in\text{Alg}_2(\mathbb{C})$ is fully dualizable if and only if it is semisimple and finitely generated over $\mathbb{C}$ as a module. \end{thm}
\bpf
The evaluation map is $A:A\otimes A^o\rightarrow k$. This is equivalent to asking that the functor $F = {_k}A_{A\otimes A^o}\otimes (-)$ has a left and right adjoint. By theorem \ref{ConditionsForRightAdjoints}, $F$ has a right adjoint which is cocontinuous if and only if $A$ is finite dimensional over $k$. Similarly, $F$ has a left adjoint if and only if $A$ is finitely generated and projective over $A\otimes A^o$. This condition is equivalent to asking that $A$ be a finite dimensional separable algebra over $k$, and since $k$ is a perfect field, this is further equivalent to $A$ being finite dimensional and semisimple. \epf
\subsection{Calculation of the Serre Automorphism}
Mentioned earlier was how the factorization $ev^R_X=(S\otimes 1_{X^\vee})\circ coev_X$ was really a more explicit construction of the Serre automorphism. In order to obtain this factorization, it suffices to draw the coevaluation and right adjoint of evaluation in $\text{Cob}^{fr}(1)$. One will see that the map $S$ must precisely induce one twist in the normal framing of the unit interval therefore acting as the generator of $\pi_1(SO(2))$. This is carried out in \cite{schommerprieslecnotes}. Therefore, in order to write down $S$ a specific bimodule, we only need to deduce what $ev^R_X$ is. Right adjunction gives
\[ \hom_{\mathbb{C}}(ev\otimes N, M) = \hom_{\mathbb{C}}(A\otimes_{A\otimes A^o} N, M) \cong \hom_{A\otimes A^o}(N, ev^R\otimes_{\mathbb{C}} M)\]
Since these isomorphisms hold for arbitrary $N$ and $M$, pick $M=\mathbb{C}$ and $N=A\otimes A^o$. Then we have isomorphisms
\[ \hom_{\mathbb{C}}(A\otimes_{A\otimes A^o}(A\otimes A^o), \mathbb{C})\cong \hom_{A\otimes A^o}(A\otimes A^o, ev^R\otimes_{\mathbb{C}}\mathbb{C})\]
More simply, $A^\vee = \hom_{\mathbb{C}}(A,\mathbb{C}) = \hom_{A\otimes A^o}(A\otimes A^o, ev^R)=ev^R$

Now, writing out the factorization of $ev^R$ gives: $A^\vee = ev^R_X = (S\otimes 1_{X^\vee})\circ coev_X = {_A}S_A\otimes_k {_{A^o}}A^o_{A^o} \otimes_{A\otimes A^o}A$. Since $A\otimes A^o$ modules correspond to $A$-$A$ bimodules, this expression reduces to ${_A}S_A\otimes_k {_A}A_A\otimes_A {_A}A_A = S\otimes_k A$. However, this last module has the action of $A\otimes A$ with the first factor acting on $S$ and the second acting on $A$ so this tensor product is just $S$. Therefore one immediately sees that $S=A^\vee$ as an $A$-$A$ bimodule or $A\otimes A^o$-module. Since an oriented field theory is a equivalent to a fully dualizable (i.e. finitely generated semisimple) algebra $A$, the trivialization of the Serre automorphism is the same thing as giving an isomorphism of $A$-$A$ bimodules $A\rightarrow A^\vee$. This isomorphism produces a perfect pairing
\[A\otimes_{A\otimes A^o} A \rightarrow \mathbb{C}\]
Therefore, one can think of this trivialization of the Serre automorphism $S$ as giving a trace $tr:A\otimes_{A\otimes A^o} A\rightarrow \mathbb{C}$. To summarize, we have the following classification result:
\begin{thm}
\label{FullStatementOfClassification}
\emph{(Classification of $2$ dimensional Extended Topological Field Theories)}
An extended $2$ dimensional framed topological field theory $Z$ is specified up to natural equivalence by a finitely generated semisimple algebra $A$ over $\mathbb{C}$ up equivalence in $\mathrm{Alg}_2(\mathbb{C})$ (this equivalence is called ``Morita equivalence''). The data of an oriented $2$ dimensional topological field theory is just providing a nondegenerate bilinear form on $A$ so that $2$ dimensional oriented topological field theories are equivalent to the data of a finitely dimensional semisimple Frobenius algebra $A$.
\end{thm}
Since Morita equivalent rings have the same module categories, two Morita equivalent rings $A$ and $B$ must have the same trivializations of the Serre automorphism. Therefore, the set of traces on $A$ and $B$ turning these two rings into Frobenius algebras must biject with one another so the same oriented field theories come from Morita equivalent algebras.

Before moving on to examples of $2$ dimensional topological field theories, we take a quick detour to explain the nomenclature of the Serre automorphism. It turns out that in general, it is important to study when a $k$-linear category $\mathcal{C}$ has a functor $S:\mathcal{C}\rightarrow \mathcal{C}$ such that $\hom(X,SY)=\hom(Y,X)^\vee$. If this functor exists, it is unique up to canonical isomorphism. Furthermore, if we restrict our attention only to the subcategory generated by ``compact objects'' of $\mathcal{C}$, the existence of the Serre automorphism implies the full dualizability of $\mathcal{C}$ in the $2$-category of linear categories, functors and natural transformations \cite{teleman12}. In the case where $\mathcal{C}=\text{D}^+\text{QCoh}(X)$ the lower-bounded derived category of quasicoherent $\mathcal{O}_X$-modules for a smooth projective variety $X$ of dimension $d$, Serre Duality allows us to explicitly construct the Serre automorphism.
\begin{thm}
\label{SerreDuality}
\emph{(Serre Duality)}
Let $\mathcal{F}\in\mathcal{C}$ be a quasicoherent $\mathcal{O}_X$-module, then there is a perfect pairing
\[ \hom_{\mathcal{C}}(\mathcal{F}, \omega_X[d]) \times \hom_{\mathcal{C}}(\mathcal{O}_X,\mathcal{F})\rightarrow k \]
inducing an isomorphism $\hom_{\mathcal{C}}(\mathcal{F}, \omega_X[d]) \cong \hom_{\mathcal{C}}(\mathcal{O}_X,\mathcal{F})^\vee$.\end{thm}
Now, set the functor $S$ so that $S\mathcal{F}=\mathcal{F}\otimes\omega_X[d]$. Then we have
\begin{eqnarray}
\hom(\mathcal{F},S\mathcal{G}) && = \hom(\mathcal{F},\mathcal{G}\otimes \omega_X[d])\nonumber \\
&& =\hom(\mathcal{F}\otimes \mathcal{G}^\vee , \omega_X[d])\nonumber \\
&& =\hom(\mathcal{O}_X,\mathcal{F}\otimes\mathcal{G}^\vee)^\vee \nonumber \\
&& =\hom(\mathcal{G},\mathcal{F})^\vee \nonumber
\end{eqnarray}
Therefore, we see that $S$ must be the Serre functor and if we take only the complexes in $\mathcal{C}$ which are compactly generated, then $\mathcal{C}$ is a fully dualizable object. The reader  interested in what framed field theory this category defines is referred to \cite{teleman12}; details of Serre duality can be found in almost any standard text on algebraic geometry, e.g. Hartshorne \cite{hartshorneres}\cite{hartshorneAG}.
\section{Extended Dijkgraaf-Witten Theory (Finite Gauge Theory)}
\subsection{Construction}
We can now attempt to build the corresponding extended topological quantum field theory for group rings $\mathbb{C}[G]$ for a finite group $G$. These rings are all semisimple since they decompose into direct sums of the irreducible representations of $G$ and so by theorem \ref{FullyDualizableAlgebras}, these should correspond to extended $2$ dimensional field theories. While at first, it seems that group rings might be very distinguished examples of semisimple algebras, we will show that the theories arising from these group rings, so called Dijkgraaf-Witten theories or finite gauge theories, account for every possible extended $2$ dimensional field theory up to natural equivalence of functors.
\begin{df} For a finite group $G$, the extended topological quantum field theory which evaluates to $\mathbb{C}[G]$ for the positively oriented point is called ``extended Dijkgraaf-Witten theory'' or extended finite gauge theory for $G$.\end{df}
To begin, we first construct the framed version of Dijkgraaf-Witten theory. Then, we know from the Cobordism Hypothesis with group structure that to get an oriented theory only requires the additional specification of a trace map on $\mathbb{C}[G]$ (different choices of trace maps can turn the underlying framed theory into different oriented theories). By the Cobordism Hypothesis, we know that we must prescribe $A=Z(+)=\mathbb{C}[G]$. Now, consider the evaluation map $\emptyset\rightarrow +\coprod -$ and the coevaluation map $+\coprod -\rightarrow \emptyset$ (the category is symmetric monoidal so we may freely swtich the ordering on disjoint unions). The framing of the evaluation and coevaluation map is just the standard framing on the unit interval now bent into a half circle (so that the composition has the framing on $S^1$ that precisely is the one not inherited from $\mathbb{R}^2$). Since the evaluation and coevaluation map are $A_{A\otimes A^o}$ and $_{A\otimes A^o}A$, we get that $Z(S^1)=A \otimes_{A\otimes A^o} A$. Now, suppose $f\in A$. Consider the element $f\otimes 1$ in $A \otimes_{A\otimes A^o} A$. Since $1=gg^{-1}$, it must be that $f=g^{-1}fg$ for any $g\in G$. Therefore, $f$ has to be of the form $\sum a_i H_i$ where $a_i\in\mathbb{C}$ and $H_i$ is the sum of all the elements in a conjugacy class of $G$. One can think of this as assigning a complex number $a_i$ to each conjugacy class $H_i$ and we recover that $Z(S^1)$ is the vector space of class functions from $G$ to $\mathbb{C}$. At this moment, one may think that this identification also automatically gives a counit (and therefore trace) map: intuitively, $Z(S^1)$ has a canonical identity element which is just the function assigning the number $1$ to the conjugacy class of the identity. However, this map does not correspond to a framed bordism in $\text{Cob}(2)$ since the framing of the circle not coming from $\mathbb{R}^2$ does not extend over the disc.

The more difficult construction is to determine what morphism of modules to assign to $2$-manifolds. Since the pair of pants is the fundamental building block of all other $2$-bordisms, it suffices to check to see what morphism it gives from $Z(B): Z(S^1\coprod S^1)\rightarrow S^1$. The method is as follows: a conjugacy class of $G$ is really just a conjugacy class of homomorphisms $\mathbb{Z}\rightarrow G$. Since $\mathbb{Z}=\pi_1(S^1)$, this is the same thing as specifying a $G$-bundle on $S^1$. Therefore, a class function is just a way of assigning a complex number to every isomorphism class of $G$-bundles on $S^1$. Now, consider the following diagram corresponding to restriction of $G$-bundles:
\[ \xymatrix{
			&Maps(B,BG)\ar[rd]^{p_+}\ar[ld]_{p_-} \\
Maps(S^1, BG) && Maps(S^1\coprod S^1, BG)
}\]
Each connected component of each mapping space gives an isomorphism class of a principal $G$-bundle on the source of the mapping space so the vector spaces $Z(S^1)$ and $Z(S^1\coprod S^1)$ are really the vector spaces of locally constant complex-valued functions on each of these mapping spaces. Then given a locally constant function $f:Maps(S^1\coprod S^1,BG)\rightarrow \mathbb{C}$, we can pull this function back to $Maps(B,BG)$ as follows: given an isomorphism class of $G$-bundle $\xi$ on $B$, we can restrict it to the boundary component $S^1\coprod S^1$ and get a bundle $\xi|_{S^1\coprod S^1}$. Then we define a function $g(\xi) = f(\xi|_{S^1\coprod S^1})$. Since $\pi_1(B)$ for the pair of pants $B$ is the free group on two generators $\mathbb{Z}\ast\mathbb{Z}$, a principal $G$-bundle $\xi$ is the same as a conjugacy class of a pair $(g,h)\in G\times G$. Indeed, projecting into the first and second factors produces conjugacy classes giving the restricted bundle $\xi|_{S^1\coprod S^1}$ for which we have a way of assigning complex values. Therefore, the only nontrivial map to define is the map $p_-$: given a rule for assigning complex numbers to $G$-bundles $\xi$ on $B$, one must find a way of assigning complex numbers to bundles on the incoming boundary $S^1$. The algorithm is as follows: given a connected component of $Maps(S^1, BG)$, look at the preimage in $Maps(B,BG)$. For each connected component $U_i$, we get a $G$-bundle $\xi_i$. Then we compute the number
\[a = \sum_{i} \frac{g(\xi_i)}{\#\pi_1(U_i)}\]
And assign the restricted bundle on $S^1$ the number $a$. In practice, this computation is actually not as difficult as it might look: the first thing to do is to label the generators of $\pi_1(B)$. Label one of the generators $\gamma$ and the other $\gamma'$; these correspond to a single loop around the outgoing boundary $S^1\coprod S^1$ of the bordism. Then the loop around the incoming boundary $S^1$ is just $\gamma'\gamma$. Therefore, given a $G$-bundle $\zeta\in Maps(S^1, BG)$ or equivalently a conjugacy class $F\subseteq G$, we look at each way we can write $F$ as the product of two conjugacy classes $F=HH'$. Then we add together the value $g(H,H')$ for all the pairs $H$ and $H'$ multiplying to $F$.

Recall that the way we were able to extract class functions from elements $A\otimes_{A\otimes A^o} A$ was to write elements as the form $\sum a_iH_i$ and set $f(H_i)=a_i$. With this in mind, the procedure we just described is exactly the way polynomials are multiplied: in particular, given two class functions $f=\sum a_iH_i$ and $g=\sum b_iH_i$, the image in $Z(S^1)$ under $Z(B)$ is $fg=\sum a_iH_i\sum b_iH_i$ which is exactly the multiplication inherited from the group ring $\mathbb{C}[G]$.

If $B$ is a manifold without boundary, the formula reduces to
\[Z(B)=\sum_{U_i\in\pi_0 Maps(B,BG)} \frac{1}{\#\pi_1(U_i)}\]
Since each connected component $U_i$ gives a principal $G$-bundle $\xi_i$, this formula is the same thing as counting each isomorphism class of principal $G$-bundle weighted by the order of its automorphism group $\pi_1(U_i)=G$.
Now, we can ask what kinds of manifolds our framed theory allows us to compute invariants for. Unfortunately, the only $2$-manifold that has a trivializable tangent bundle is the torus $S^1\times S^1$. Therefore, it is desirable to give this theory additional structure in order to make it an oriented field theory. Such additional structure amounts to a trace map $A\otimes_{A\otimes A^o} A\rightarrow \mathbb{C}$ which induces a perfect pairing therefore giving an isomorphsim of $A\otimes A^o$ modules $A\rightarrow A^\vee$. Therefore, we can look at the central elements $z(A)=\hom_{A\otimes A^o}(A,A)=\hom(k,A^\vee\otimes A)$ and so we are given an isomorphism $z(A)\rightarrow A\otimes_{A\otimes A^o} A$ albeit non canonically. As a result, we see that $z(A)$ has a non-degenerate trace form and so is a commutative Frobenius algebra and hence gives an ordinary topological quantum field theory on oriented manifolds of dimensions $1$ and $2$.

\subsection{The Morita Theory of Semisimple Algebras}
Since there is an equivalence of categories between the category of framed topological quantum field theories with values in algebras bimodules and intertwiners to the category of fully dualizable objects in $\text{Alg}_2(\mathbb{C})$, it suffices to check if two fully dualizable objects are equivalent in order to check if the field theories they define are naturally isomorphic. Equivalence in the category $\text{Alg}_2(\mathbb{C})$ is Morita equivalence.
\begin{df} Two algebras $A$ and $B$ are said to be ``Morita equivalent'' if the categories of right modules are equivalent under a cocontinuous functor. More explicitly, the two algebras are Morita equivalent if there is a $B$-$A$ bimodule ${_B}P_A$ such that the functor ${_B}P_A\otimes (-)$ is an equivalence of categories. \end{df}
Note that Morita equivalences will have both adjoints since modules ${_B}P_A\otimes(-)$ inducing Morita equivalences will be finitely presented and projective \cite{lam2}.

It turns out the Morita theory for semisimple algebras is actually quite simple: this is because of a powerful structure theorem for semisimple algebras.
\begin{thm}
\label{ArtinWedderburn}
\emph{(The Artin-Wedderburn Theorem)}
Let $A$ be a semisimple algebra. Then $A$ is isomorphic to a product of matrix rings over division rings:
\[ A\cong \prod M_{n_i}(D_i) \]
This decomposition is unique up to order of the factors \cite{lam}.
\end{thm}
Since our algebras $A$ are all over $\mathbb{C}$, the decomposition reduces to giving $A$ as the product of matrix rings with entries in $\mathbb{C}$. This decomposition combined with a standard result from Morita theory allows us to fully classify Morita equivalences of semisimple algebras.
\begin{thm} For a ring $R$, any matrix ring $M_m(R)$ is Morita equivalent to $R$ \cite{lam2}.\end{thm}
These two theorems show that any semisimple algebra $A$ is Morita equivalent to a ring that is some product of the scalar field $\mathbb{C}$. In particular, since each copy of the scalar field embeds to the matrix ring $M_m(\mathbb{C})$ by multiples of the identity matrix, this ring $\mathbb{C}^n$ is actually the center $z(A)$. Therefore, any semisimple algebra $A$ is Morita equivalent to $z(A)$. Since $z(A)$ is commutative, Morita equivalence is the same as isomorphism, so two semisimple rings are Morita equivalent if and only if they have isomorphic centers.
\begin{cor} Let $G$ be a finite group and $n$ be the number of conjugacy classes of $G$. Then $\mathbb{C}[G]$ is Morita equivalent to $\mathbb{C}[\mathbb{Z}/n]$.\end{cor}
\bpf
The group ring $\mathbb{C}[G]$ is essentially the regular representation for $G$. The regular representation is the sum of all of the irreducible representations of $G$ with each irreducible representation contributing a factor $M_{n_i}(\mathbb{C})$ in the Artin-Wedderburn decomposition of $\mathbb{C}[G]$. Since the number $n$ of irreducible representations of a finite group $G$ is the same as the number of its conjugacy classes, the center of $\mathbb{C}[G]$ is $\mathbb{C}^n$. Since $\mathbb{Z}/n$ is Abelian, it clearly also has $n$ conjugacy classes and therefore the center of its group ring is also $\mathbb{C}^n$ and we see that $\mathbb{C}[G]$ is Morita equivalent to $\mathbb{C}[\mathbb{Z}/n]$. \epf
\begin{cor}A $2$ dimensional framed extended topological quantum field theory with values in algebras bimodules and intertwiners are naturally isomorphic to finite gauge theory for a cyclic group. \end{cor}
\section{Conclusion}

Putting together all the results we have so far, we can provide a concise description of all $2$ dimensional topological field theories: suppose $Z$ is an ordinary oriented $2$ dimensional field theory, then $Z$ defines a commutative Frobenius algebra $F$. If this ring $F$ is semisimple, then it is isomorphic to a product of fields $F=\prod^n_{i=1} k$. Any two semisimple algebras $A$ and $A'$ over $k$ that has $F$ as their center will be Morita equivalent: this is because $A$ and $A'$ must have the same number of factors in their Artin-Wedderburn decomposition since each factor contributes a copy of $k$ to the center. Therefore there will be a unique extended $2$ dimensional framed field theory extending $Z$ given by $A$ and we might as well pick $A=k[\mathbb{Z}/n]$. The possible oriented field theories coming from this underlying framed field theory are in bijection with $A$-$A$ bimodule isomorphisms $A^\vee\rightarrow A$ and each trace gives an isomorphism of $Z(S^1)=A\otimes_{A\otimes A^o}A$ with the center of $A$. These results allow one to really understand all topological field theories in $2$ dimensions: we have an easily computable criterion for when the field theory has an extension and we can list and explicitly construct every fully extended $2$ dimensional field theory. This is largely due to the fact that full dualizability in $2$-categories do not require us to write down too many conditions; as expected, in $3$ dimensions the story is different. Instead of having a semisimple algebra give an extended field theory, one has the much more complicated ``fusion category'' \cite{kirillov}. Dijkgraaf-Witten theory happens to translate nicely over to the $3$ dimensional case to produce fully extended field theories, but in general, there are problems in extending other well known field theories such as Chern-Simons theory \cite{lht}. These $3$ dimensional theories give rise to more interesting phenomena than their $2$ dimensional counterparts, for example, Witten has shown that knot invariants are related to $3$ dimensional field theories in \cite{witten89}; another common application of $3$ dimensional topological field theories is build two equivalent theories from different constructions and compare the invariants produced. For $(\infty,n)$-categories with $n$ higher than $4$, not much is known about full dualizability and the subject remains largely uncharted.

\clearpage
\renewcommand{\abstractname}{Acknowledgements}
\begin{abstract}
I would like to thank my adviser Jacob Lurie without whom this thesis would have been impossible. He provided invaluable advice and insights and patiently explained many concepts and calculations to me during our weekly meetings throughout this research project. I would also like to acknowledge Zheng-Han Wang at Microsoft Station Q and Xiao-Gang Wen at MIT for inviting me to a summer workshop at Peking University on topological quantum computing where I was first exposed to topological field theories. I would also like to thank Dung-Hai Lee for our many conversations clarifying the physical intuition and motivation for topological field theories. Finally, I thank Constantin Teleman and Christopher Schommer-Pries for our email exchanges during my writing.
\end{abstract}

\end{document}